\documentclass[11pt,a4paper,reqno]{amsart}
\usepackage{amsfonts}
\usepackage{amsmath,amssymb,amsthm,amsxtra}
\usepackage{float}
\usepackage[colorlinks, linkcolor=blue,anchorcolor=Periwinkle,
    citecolor=red,urlcolor=Emerald]{hyperref}
\usepackage[usenames,dvipsnames]{xcolor}
\usepackage{enumitem}
\usepackage{geometry,array} 
\usepackage{graphicx}
\usepackage{subfigure}
\usepackage{bookmark}
\usepackage{tikz}\usetikzlibrary{matrix}
\usepackage{url}

\usetikzlibrary{decorations.markings}
\tikzset{->-/.style={decoration={  markings,  mark=at position #1 with
    {\arrow{>}}},postaction={decorate}}}
\tikzset{-<-/.style={decoration={  markings,  mark=at position #1 with
    {\arrow{<}}},postaction={decorate}}}
\usepackage[all]{xy}
\usetikzlibrary{arrows}
\usetikzlibrary{decorations.pathreplacing,decorations.pathmorphing,shadings,fadings,calc}
\usepackage{extarrows}
\newcommand{\arxiv}[1]{\href{http://arxiv.org/abs/#1}{arXiv:#1}}

\theoremstyle{plain}
\newtheorem{theorem}{Theorem}[section]

\newtheorem*{con/thm}{Construction/Theorem}

\newtheorem{proposition}[theorem]{Proposition}
\newtheorem{conjecture}[theorem]{Conjecture}
\theoremstyle{definition}
\newtheorem{definition}[theorem]{Definition}

\newtheorem{example}[theorem]{Example}

\newtheorem{remark}[theorem]{Remark}

\numberwithin{equation}{section}

\def\hua{\mathcal}
\def\kong{\mathbb}
\def\<{\langle}
\def\>{\rangle}

\def\Aut{\operatorname{Aut}}

\def\Hom{\operatorname{Hom}}

\def\Ext{\operatorname{Ext}}

\def\Stab{\operatorname{Stab}}
\def\Stap{\operatorname{Stab}^\circ}

\newcommand{\h}{\operatorname{\hua{H}}}            

\renewcommand{\k}{\mathbf{k}}

\newcommand{\tilt}[3]{{#1}^{#2}_{#3}}
\newcommand{\Cone}{\operatorname{Cone}}

\def\numbers{\begin{enumerate}[label=\arabic*{$^\circ$}.]}
\def\ends{\end{enumerate}}

\newcommand{\EG}{\operatorname{EG}}       
\newcommand{\SEG}{\operatorname{SEG}}       
\newcommand{\C}{\hua{C}}

\newcommand{\CEG}[2]{\operatorname{CEG}_{#1}(#2)}             
\newcommand{\D}{\operatorname{\hua{D}}}
\def\eg{\operatorname{\hua{EG}}}

\def\ceg{\operatorname{\hua{CEG}}}
\def\conj{\operatorname{ad}}

\newcommand{\per}{\operatorname{per}}
\def\zero{\hua{H}_\Gamma}
\newcommand{\Tri}{\bigtriangleup}

\def\arrow{red}
\def\surf{\mathbf{S}}                       

\newcommand{\DT}{\operatorname{DT}}        
\newcommand{\ST}{\operatorname{ST}}        
\newcommand{\BT}{\operatorname{BT}}        

\newcommand{\MCG}{\operatorname{MCG}}

\def\coh{\operatorname{Coh}}

\def\TT{\mathbf{T}}
\def\T{\kong{T}}

\def\M{\mathbf{M}}

\newcommand{\AT}{\operatorname{AT}}        

\def\surfo{{\mathbf{S}}_\Tri}

\def\cA{\operatorname{CA}}

\newcommand\Bt[1]{\operatorname{B}_{#1}}

\newcommand{\Quad}{\operatorname{Quad}}

\newcommand{\SBG}{\operatorname{SBr}}

\newcommand\Br{\operatorname{Br}}
\newcommand\CBr{\operatorname{CT}}
\newcommand\Co{\operatorname{Co}}

\newcommand\Tr{\operatorname{Tr}}

\newcommand{\qq}[1]{\operatorname{\Gamma}_{#1}Q}
\usepackage{braids}

\title{The braid group for a quiver with superpotential}
\author[Qiu]{Yu Qiu}
\address{YQ:
Department of Mathematics,
Chinese University of Hong Kong,
Shatin, N.T.,Hong Kong}
\email{yu.qiu@bath.edu}

\date{\today}
\begin{document}

\begin{abstract}
We survey various generalizations of braid groups for quivers with superpotential
and focus on the cluster braid groups, which are introduced in a joint work with A.~King.
Our motivations come from the study of cluster algebras,
Calabi-Yau categories and Bridgeland stability conditions.

\end{abstract}


\maketitle


\def\dgen{2\mathbb{N}_{\leq g}-1}
\def\dgene{2\mathbb{N}_{\leq g}}

\section{Introduction}
\subsection{Cluster algebras}
Cluster algebras were introduced by Fomin-Zelevinsky \cite{FZ}
whose original motivation comes from the study of
total positivity in algebraic groups and canonical bases in quantum groups.
For the last two decades, the theory of cluster algebras has grown exponentially
due to such a phenomenon appears in many subjects, such as
Poisson geometry, integrable systems, Teichm\"{u}ller spaces,
algebraic geometry, mirror symmetry, representation theory of algebras...
(cf. Keller's survey \cite{K0}).

The combinatorial aspect of the cluster theory is quiver mutation,
which was further developed by Derksen-Weyman-Zelevinsky \cite{DWZ}
as mutation of quivers with potential.
A cluster algebra will determine a class of mutation-equivalent quivers of potential.
Our aim is to explain what is the (generalized) braid group that should be associated to
a quiver with potential.
There are various answers depending on motivations.
Nevertheless, our motivation comes from
the study of Bridgeland stability conditions on Calabi-Yau categories that are associated to
cluster algebras/quivers with potentials.
Therefore, our criterion of introducing such braid groups
is to serve the study of the topology of spaces of stability conditions.

\subsection{Stability conditions on Calabi-Yau categories}
Stability conditions on triangulated categories were introduced by Bridgeland \cite{B1},
which were motivated from Douglas' $\Pi$-stability in the study of D-branes in string theory.
The crucial feature is that the set of all stability conditions on a triangulated category $\D$
is in fact a complex manifold $\Stab\D$.
Many interesting examples of $\D$ are from geometry,
namely those appear in homological mirror symmetry
\begin{gather}\label{eq:HMS}
    \operatorname{DFuk}(\mathbf{X})\cong\D^b(\coh\mathbf{Y}),
\end{gather}
where $\mathbf{X}$ is a symplectic mainfold on A-side and $\mathbf{Y}$ its complex mirror on B-side
(cf. \cite{KS,ST}).

While in general the study of the space $\Stab\D$ on the Calabi-Yau-3 categories
mentioned above is very hard, there are many progress in the simplified quivery cases.
Namely, for a quiver with potential $(Q,W)$ (say from cluster algebra setting),
one can construct a Calabi-Yau-3 category $\D_{fd}(\Gamma(Q,W))$ via Ginzburg dg algebra $\Gamma(Q,W)$.
In some of the setting, e.g.
quivers with potential 
from a marked surface $\surf$ in the sense of Fomin-Shapiro-Thurston
and Labardini (\cite{FST,LF}),
such categories (only depend on $\surf$)
\[
    \D(\surf)\colon=\D_{fd}(\Gamma(Q,W))
\]
can be embedded into categories in \eqref{eq:HMS}
(due to Smith \cite{S}).
Moreover, Bridgeland-Smith \cite{BS} establish a connection between
Teichm\"{u}ller theory and stability conditions.
More precisely, they prove that
\begin{gather}\label{eq:motivation}
    \Stap\D(\surf)/\Aut\D(\surf)\cong\Quad(\surf),
\end{gather}
where $\Quad(\surf)$ is the moduli space of (signed) quadratic differentials on $\surf$.
Their motivations are coming from string theory in physics,
Donaldson-Thomas theory and (homological) mirror symmetry (cf. \cite{GMN}, \cite{S} and \cite{QQ}).

To the symmetry groups in the formula \eqref{eq:motivation},
one needs to understand the spherical twist group
$\ST\D(\surf)\subset\Aut\D(\surf)$ that sits in the short exact sequence
\[
    1\to\ST\D(\surf)\to\Aut\D(\surf)\to\MCG(\surf)\to1,
\]
where $\MCG(\surf)$ is the mapping class group of $\surf$.
Such spherical twist groups were first study by Khovanov-Seidel-Thomas \cite{KS,ST}
from the two sides of the homological mirror symmetry
in the case when $\surf$ is a disk.
In the previous works on spherical twist groups (e.g. \cite{KS,ST,BT,QW}),
one usually proved that $\ST\D(\surf)$ is isomorphic to the braid group of
the corresponding (Dynkin) type.
However, in the general case (arbitrary surfaces or arbitrary quivers with potential from cluster algebras), the associated braid groups are not (well-)defined yet.

\subsection{Context}
In Section~\ref{sec:2}, we introduce the classical braid groups and several generalizations
from different point of views.
In particular, we will review the symplectic generalization via
spherical twists on Calabi-Yau categories in Section~\ref{sec:ST-CY}.
We summarize the previous results in Section~\ref{sec:pre}.
In Section~\ref{sec:4}, we introduce the cluster braid groups, which is due to the forthcoming joint work with Alastair King \cite{KQ2}.
Such a generalization, via cluster exchange groupoid, can apply to higher Calabi-Yau cases,
i.e. in the quivers with superpotential setting.
Further study in this direction will appear in the project joint with Akishi Ikeda and Yu Zhou.

Some conventions:
\begin{itemize}
\item The convention of composition is from left to right (as product).
\item $\Co(a,b)\Longleftrightarrow ab=ba$.
\item $\Br(a,b)\Longleftrightarrow aba=bab$.
\item $\Tr(a,b,c)\colon abca=bcab=cabc$.
\end{itemize}

\subsection*{Acknowledgments}
I'd like to thank my collaborators Alastair King, Jon Woolf and Yu Zhou
for collaborating with me on the topic of braid groups.
I also want to thank Fang Li, Zongzhu Lin and Bin Zhu,
the organizers of the International Workshop on Cluster Algebras in
Naikai University, Tianjin (2017),
for inviting me writing this article contributing to the conference proceedings.

\section{Various generalization of braid groups}\label{sec:2}
\subsection{An algebraic generalization}
The \emph{classical braid group} (a.k.a. Artin group) $\Br_{n+1}$
on $n+1$ strands has the following presentation
\[\begin{array}{rll}
    \Br_{n+1}=\< b_1, \ldots, b_n \mid &\Br(b_i,b_{i+1})\;
        \forall 1\leq i\leq n-1,\;\\& \Co(b_j,b_k) \;\forall\;|j-k|>1 \>.
\end{array}\]
\begin{definition}\label{def:Br}
The \emph{(generalized) braid group} $\Br(Q)$ associated to a quiver $Q$
(or its underlying diagram $\underline{Q}$) is
the group with generators $b_i, i\in Q_0$, and the relations
\[\begin{cases}
    \Co(b_i,b_j), &\text{if there are no arrows between $i$ and $j$},\\
    \Br(b_i, b_j), &\text{if there is exactly an arrow between $i$ and $j$}.
\end{cases}\]
When $Q$ is of type $A_n$, we have $\Br(Q)=\Br_{n+1}$.
\end{definition}

A potential $W$ of a quiver is the sum of certain cycles in $Q$.
According to the philosophy in \cite{QQ},
the (proper) braid group $\Br(Q,W)$ of a quiver with potential $(Q,W)$
should admit the following presentation:
\begin{itemize}
\item generators are (indexed by) $Q_0$;
\item commutation/braid relations correspond to zero/exactly one arrow between vertices.
\item a triangle relation $\Tr$ for each 3-cycle $abc$ in the potential term.
\end{itemize}
Such a definition is not `correct' in general; however, at least we have the following.

\begin{definition}\label{def:AT}
Suppose that $(Q,W)$ is a quiver with potential
such that there is at most one arrow between any two vertices.
The algebraic braid twist group $\AT(Q,W)$ is defined by the presentation
\begin{itemize}
\item generators $b_i$ are (indexed by) $i\in Q_0$;
\item there is a relation $\Br(b_i,b_j)$ if there is an arrow between $i,j$;
otherwise there is a relation $\Co(b_i,b_j)$.
\item there are relations $R_i=R_j$ for any $i\neq j$ ,
  if there is a cycle $Y\colon1 \to 2 \to \cdots \to m \to 1$ in $W$,
  where $R_i=b_i b_{i+1} \cdots b_{2m+i-3}$ with convention $k=m+k$ here.
  Note that for any $i\neq j$ and $i'\neq j'$, we have $R_i=R_j\Longleftrightarrow R_i'=R_j'$.
\end{itemize}
\end{definition}

\subsection{A geometric generalization}
\begin{figure}[h]
\begin{tikzpicture}[xscale=2,yscale=1]
\draw[dashed, fill=gray!11] (3,0) ellipse (2.3 and .8);
\draw[dashed, fill=gray!11] (3,-9.5) ellipse (2.3 and .8);
\braid[thick, style strands={1}{ForestGreen},
style strands={2}{blue},
style strands={3}{violet},
style strands={4}{orange},
style strands={5}{green}]
s_1 s_3 s_4 s_2^{-1} s_1 s_3 s_2^{-1} s_1 s_2^{-1};
\foreach \j in {3,4,5,2,1}
{\draw (\j,0)node[white]{$\bullet$}node[red]{$\circ$}
    (\j,-9.5)node[white]{$\bullet$}node[red]{$\circ$};}
\end{tikzpicture}\caption{Classical braids}\label{fig:braids}
\end{figure}
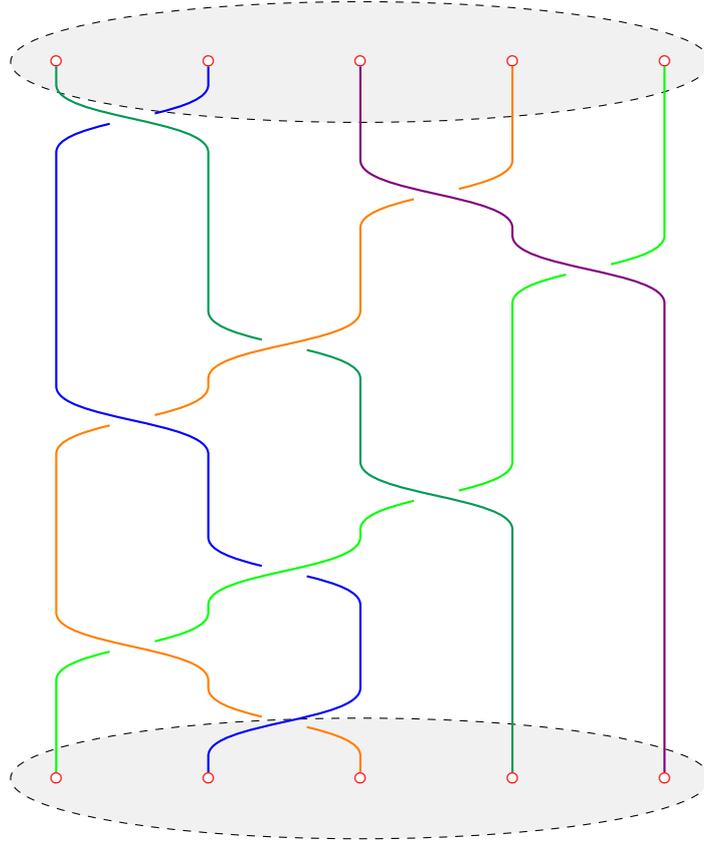

One can visualize a braid (that is, an element of the classical braid group $\Br_{n+1}$
in $\mathbb{R}^3$), as shown in Figure~\ref{fig:braids}.
The element there is
\[b_1 b_3 b_4 b_2^{-1} b_1 b_3 b_2^{-1} b_1 b_2^{-1}
    \in\Br_5.\]
Along this direction, one can define the surface braid group more precisely as follows.

Let $\surfo$ be a decorated surface, where $\surf$ is a (topological) surface with
a set $\Tri=\{Z_1,\ldots,Z_\aleph\}$ of $\aleph$ decorating points in $\surf^\circ$.
The classical case ($\Br_{n+1}$) is when $\surf$ is a disk and $\aleph=n+1$.
Then we can embed $\surfo$ in $\mathbb{R}^2$ instead of (for the general case) $\mathbb{R}^3$.

\begin{definition}\label{def:SBr}
A geometric braid on $\surfo$ based at $\Tri$ is
an $\aleph$-tuple $\Psi=(\psi_1,\ldots,\psi_\aleph)$ of paths
$$\psi_i\colon [0,1]\to\surfo$$ such that
\begin{itemize}
\item $\psi_i(0)=Z_i$;
\item $\psi_i(1)=Z_i$;
\item $\{\psi_1(t),\ldots,\psi_\aleph(t)\}$ are distinct points in $\surfo^\circ$, for $0<t<1$.
\end{itemize}
The product of geometric braids follows the same way of products of paths
(in the fundamental group setting).
All braids on $\surfo$ with the product above form the surface braid group $\SBG(\surfo)$.
\end{definition}

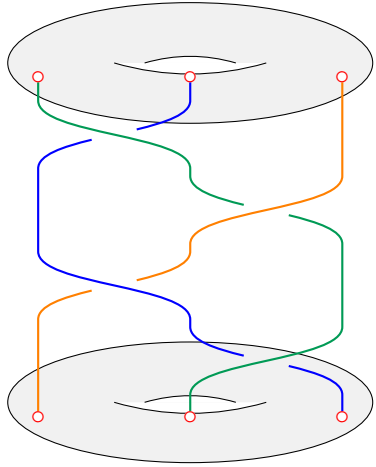
\begin{figure}[h]
\begin{tikzpicture}[xscale=2,yscale=1]
\draw[fill=gray!11] (2,0+.2) ellipse (1.2 and .8);
\draw[fill=white](2+.5,0.2)to[bend left](2-.5,0.2);
\draw[fill=white](2+.3,0.2)to[bend right](2-.3,0.2);
\draw[fill=gray!11] (2,-4.5+.2) ellipse (1.2 and .8);
\draw[fill=white](2+.5,0.2-4.5)to[bend left](2-.5,0.2-4.5);
\draw[fill=white](2+.3,0.2-4.5)to[bend right](2-.3,0.2-4.5);
\braid[thick, style strands={1}{ForestGreen},
style strands={2}{blue},
style strands={3}{orange}]
s_1 s_2^{-1} s_1 s_2^{-1};
\foreach \j in {3,2,1}
{\draw (\j,0)node[white]{$\bullet$}node[red]{$\circ$}
    (\j,-4.5)node[white]{$\bullet$}node[red]{$\circ$};}
\end{tikzpicture}\caption{Torus braids (cf. \cite{GJP})}\label{fig:braid2}
\end{figure}

For instance, when $\surfo$ is a torus with three decorations,
Figure~\ref{fig:braid2} tries to show a braid on $\surfo$.
\subsection{A topological generalization}\label{sec:top}
Another well-known alternative definition of the classical braid group $\Br_{n+1}$ is the following:
\[
    \Br_{n+1}=\pi_0\operatorname{Diff}(\mathbf{D}_{n+1})=\colon\operatorname{MCG}(\mathbf{D}_{n+1}),
\]
where $\mathbf{D}_{n+1}$ is a close disk with a set $\Delta$ of $n+1$ decorations (points)
and  $\operatorname{Diff}(\mathbf{X})$ is the
the group of diffeomorphisms of $\mathbf{X}$ that preserve
the boundary pointwise and $\Delta$ setwise.

The corresponding generalization can be done in the following way.
A closed arc in $\surfo$ is a curve (up to isotopy) in $\surf$ whose interior lies in $\surf-\Tri$ and whose endpoints are different decorating points in $\Tri$.
Denote by $\cA(\surfo)$ the set of simple closed arcs in $\surfo$.
For any closed arc $\eta\in\cA(\surfo)$, there is the (positive) braid twist
$\Bt{\eta}\in\MCG(\surfo)$ along $\eta$, which is defined by Figure~\ref{fig:Braid twist}.
\begin{definition}
The braid twist group $\BT(\surfo)$ of $\surfo$ is the subgroup of $\MCG(\surfo)$
generated by $\{\Bt{\eta}\mid\eta\in\cA(\surfo)\}$.
\begin{figure}[ht]\centering
\begin{tikzpicture}[scale=.2]
  \draw[very thick,NavyBlue](0,0)circle(6)node[above,black]{$_\eta$};
  \draw(-120:5)node{+};
  \draw(-2,0)edge[red, very thick](2,0)  edge[cyan,very thick, dashed](-6,0);
  \draw(2,0)edge[cyan,very thick,dashed](6,0);
  \draw(-2,0)node[white] {$\bullet$} node[Emerald] {$\circ$};
  \draw(2,0)node[white] {$\bullet$} node[Emerald] {$\circ$};
  \draw(0:7.5)edge[very thick,->,>=latex](0:11);\draw(0:9)node[above]{$\Bt{\eta}$};
\end{tikzpicture}\;
\begin{tikzpicture}[scale=.2]
  \draw[very thick, NavyBlue](0,0)circle(6)node[above,black]{$_\eta$};
  \draw[red, very thick](-2,0)to(2,0);
  \draw[cyan,very thick, dashed](2,0).. controls +(0:2) and +(0:2) ..(0,-2.5)
    .. controls +(180:1.5) and +(0:1.5) ..(-6,0);
  \draw[cyan,very thick,dashed](-2,0).. controls +(180:2) and +(180:2) ..(0,2.5)
    .. controls +(0:1.5) and +(180:1.5) ..(6,0);
  \draw(-2,0)node[white] {$\bullet$} node[Emerald] {$\circ$};
  \draw(2,0)node[white] {$\bullet$} node[Emerald] {$\circ$};
\end{tikzpicture}
\caption{The braid twist (cf. \cite{GJP})}
\label{fig:Braid twist}
\end{figure}
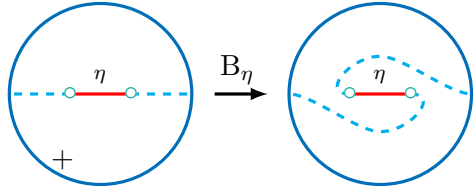
\end{definition}

Moreover, such braid twist group can be defined in the setting of quivers with potential
from marked surface (in the sense of Fomin-Shapiro-Thurston \cite{FST}, cf. \cite{LF}).
First, we equip $\surf$ with a set $\mathbf{M}$ of marked points on $\partial\surf$
such that each boundary component contains at least one marked point.
An open arc in $\surf$ is (the isotopy class of) a curve in $\surf$
that connects two marked points in $\M$, which is neither isotopic to a boundary segment nor to a point.
A triangulation is a maximal collection of open arcs such that
there are no interior intersections pairwise.
Moreover, suppose that
\[\aleph=4g_\surf+2|\partial\surf|+|\M|-4\]
and it is well-known that any triangulation $\T$ of $\surf$ contains $\aleph$ triangles.
Here $g_\surf$ is the genus of $\surf$.

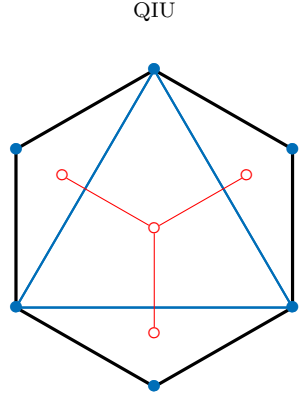
\begin{figure}[t]\centering
\begin{tikzpicture}[scale=.35,rotate=-120]
\foreach \j in {1,...,6}{\draw[very thick](60*\j+30:6)to(60*\j-30:6);}
\foreach \j in {1,...,6}{\draw[NavyBlue,thick](120*\j-30:6)node{$\bullet$}to(120*\j+90:6)
    (120*\j-90:6)node{$\bullet$};}
\foreach \j in {1,...,3}{  \draw[red](30+120*\j:4)to(0,0);}
\foreach \j in {1,...,3}{  \draw(30-120*\j:4)node[white]{$\bullet$}node[red]{$\circ$};
    \draw(30-120*\j+14:2)node[white]{\text{\footnotesize{$s_\j$}}};}
\draw(0,0)node[white]{$\bullet$}node[red]{$\circ$};
\end{tikzpicture}\caption{The dual graph of a triangulation}\label{fig:dual}
\end{figure}

Let $(Q_\T,W_\T)$ be a quiver with potential from some triangulation $\T$ of $\surf$.
Then consider a triangulation $\TT$ on $\surfo$ such that
$\TT$ becomes $\T$ when forgetting about $\Tri$ and
there is exactly one decoration in each triangle of $\T$ when adding $\Tri$ back to $\surf$.
By abuse of notation, we will not distinguish $\T$ and $\TT$ in the following.
Then \begin{itemize}
\item vertices of $Q_\TT$ are (indexed by) arcs $\gamma$ in $\TT$;
\item arrows of $Q_\TT$ correspond to angles of triangles in $\TT$;
\item terms of $W_\TT$ are 3-cycles that correspond to triangles of $\TT$.
\end{itemize}
Let $\TT^*$ be the dual graph of $\TT$, which consists of closed arcs in $\surfo$
(cf. Figure~\ref{fig:dual}).

\begin{definition}\label{def:BT}
The braid twist group $\BT(Q_\TT,W_\TT)$ associated to the quiver with potential $(Q_\TT,W_\TT)$
is the subgroups of $\MCG(\surfo)$ generated by
$\{\Bt{\eta}\mid\eta\in\TT^*\}$.
\end{definition}
If it is known that (cf. \cite{QQ,QZ3})
\[
    \BT(Q_\TT,W_\TT)\cong\BT(\surfo)\subset\SBG(\surfo).
\]

\subsection{A monodromy generalization}\label{sec:Geo}

A \emph{monodromy representation} of a braid group is
a representation of $\Br_{n+1}$, in the mapping class group of some surface, which sends the
different (standard) generators to distinct Dehn twists of certain curves on the surface.
For instance, the famous Birman-Hilden representation 
is a monodromy representation,
which is a special case (of type A) coming from the Milnor fibres.
More precisely, consider the simple singularities (with two complex variables)
\begin{eqnarray*}
  &A_n\colon& f(x,y)=x^2+y^{n+1}\qquad(n\geq1)\\
  &D_n\colon& f(x,y)=x(x^{n-2}+y^2)\qquad(n\geq4)\\
  &E_6\colon& f(x,y)=x^3+y^4\\
  &E_7\colon& f(x,y)=x(x^2+y^3)\\
  &E_8\colon& f(x,y)=x^3+y^5
\end{eqnarray*}
The Riemann surface (the singularity) $\overline{\surf}$ consists
of the singular points of the hypersurface
\[\{f(x,y)=0\}\subset\mathbb{C}^2.\]
To get the Milnor fibres (cf. \cite{PV,KS}),
one perturbs the equation $f(x,y)$, so as to smooth out the singular point,
and then intersects the outcome with a ball around the origin.
One gets a family of curves $\{C_i\mid i\in Q_0\}$
such that the intersection form between them
is given by the corresponding Dynkin diagram $\underline{Q}$ as below.
\begin{gather}\label{eq:labeling}
\begin{array}{llr}
    \xymatrix{A_{n}:\\} \quad &
    \xymatrix{
        1 \ar@{-}[r]& 2 \ar@{-}[r]& \cdots \ar@{-}[r]& n }         \\\\
    \xymatrix@R=0.1pc{\\D_{n}:\\} \quad &
    \xymatrix@R=0.4pc{
        1 \\
        & 3 \ar@{-}[ul]\ar@{-}[dl] \ar@{-}[r]& 4 \ar@{-}[r]& \cdots \ar@{-}[r]& n\\
        2\\}                                                             \\
    \xymatrix{E_{6,7,8}:} \quad &
    \xymatrix@R=1.5pc{
        && 4 \\
        1 \ar@{-}[r]& 2 \ar@{-}[r]& 3 \ar@{-}[r]\ar@{-}[u]& 5 \ar@{-}[r]& 6
                \ar@{-}[r]& 7 \ar@{-}[r]& 8}
\end{array}
\end{gather}
Then the monodromy representation $\rho_m$ of $\Br(\underline{Q})$
is given by
\begin{gather}\label{eq:rhom}
    \rho_m\colon \Br(Q)\to\pi_0\operatorname{Diff}(\overline{\surf}),
        \quad b_i\mapsto \mathrm{D}_{C_i},
\end{gather}
where $\mathrm{D}_{C_i}$ is the Dehn twist of $C_i$.
Birman-Hilden \cite{FM} proved that $\rho_m$ is faithful for type A;
Perron-Vannier \cite{PV} showed that $\rho_m$ is faithful for type D, based on
the result of Birman-Hilden.
On the other hand, Wajnryb \cite{W} showed that there is no faithful geometric representation
of the braid group of type E, which is a bit surprising.

One can generalize such an idea for the triangulated marked surfaces case.
More precisely, consider the marked surface $\surf$ mentioned above
with triangulation $\TT$.
Then the lift $\mathbf{C}_\TT$ of $\TT^*$ on the twisted surface $\Sigma_\TT$,
which is the branched double cover of $\surfo$ branching at decorations in $\Tri$,
is a collection of isotopy classes of simple closed curves
(with chosen orientation and known as clusters of curves in \cite{KQ2}).
The Dehn twist group $\DT(\mathbf{C}_\TT)$ of $\mathbf{C}_\TT$
is the subgroup of $\MCG(\Sigma_\TT)$ generated by $\{\mathrm{D}_{C}\mid C\in\mathbf{C}_\TT\}$.

Note that $\Sigma_\TT$ only depends on $\surfo$ (topologically);
however its combinatorial construction in \cite{KQ2} depends on $\TT$.
Consider the punctured case of $\surf$, i.e. adding a set $\mathbf{P}$ of punctures
on $\surf$ (which serves different purpose than decorations)
In fact, punctures/decorations on $\surf$ are poles (of order two)/zeroes (of order one)
of quadratic differentials on the Riemann surface associated to $\surf$
(further details concerning quadratic differentials cf. \cite{BS,KQ2}).
We have the following.

\begin{definition}\label{def:int.quiver}
A cluster of curves $\mathbf{C}$ is a collection of (isotopy classes of) oriented curves
(on some surface).
The (geometric) intersection quiver $Q(\mathbf{C})$ of $\mathbf{C}$
is the quiver whose vertices are curves in $\mathbf{C}$ and whose edges are
bijective to the positive geometric intersections between curves in $\mathbf{C}$.
\end{definition}
\begin{theorem}\cite{KQ3}
Given a quiver with potential $(Q_\TT,W_\TT)$ from a (tagged) triangulation $\TT$ of $\surf$,
there exists a cluster of curves $\mathbf{C}_\TT$ on certain surface $\Sigma_\TT$
(twisted surface in \cite{KQ3})
such that $Q_\TT=Q(\mathbf{C}_\TT)$.
\end{theorem}
\begin{remark}
The terms (cycles in the quiver $Q_\TT$) in $W_\TT$ correspond to `contractible polygons'
formed by curves in $\mathbf{C}_\TT$ (cf. Figure~\ref{fig:polygon}).
However, such a statement only makes sense when choosing representatives
in the isotopy classes of curves.
\begin{figure}[b]\centering
\begin{tikzpicture}[scale=.8]
\draw[white ] (0,0) circle (3);
\foreach \j in {1,...,5}{
    \draw[red,thick] (90+\j*72-8:3)edge[bend left=40,->,>=latex](90+\j*72+72+8:3);}
\end{tikzpicture}
\caption{A contractible polygon}
\label{fig:polygon}
\end{figure}
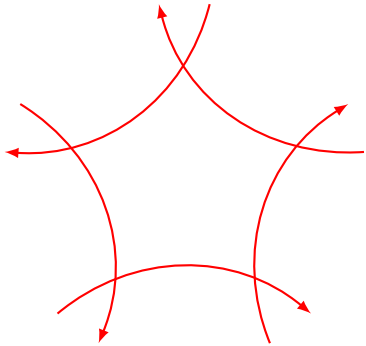
\end{remark}

\begin{definition}\cite{KQ2}\label{def:DT}
The Dehn twist group $\DT(Q_\TT,W_\TT)$ of $(Q_\TT,W_\TT)$,
is the subgroup of $\MCG(\Sigma_\TT)$ generated by $\{\mathrm{D}_{C}\mid C\in\mathbf{C}_\TT\}$.
\end{definition}

Note that in the unpunctured case, the twisted surface $\Sigma_\TT$
is exactly the branched double cover of $\surfo$, branching at $\Tri$.
Then the clusters of curves consists of curves, which are lifts of closed arcs $\TT^*$.
Then by the famous result of Birman-Hilden (\cite{FM}),  
$\DT(\Sigma_\TT)\cong\BT(Q_\TT,W_\TT)$.

\subsection{A symplectic generalization}
The symplectic representation of the braid groups arose
in the study of Kontsevich's (homological) mirror symmetry.
On the symplecitc geometry side,
Khovanov-Seidel \cite{KS} studied a subcategory $\mathcal{D}(\qq{N})$ of the derived Fukaya category of
the Milnor fibre of a simple singularities of type A.
They showed that there is a faithful braid group action on $\mathcal{D}(\qq{N})$,
where the braid group is generated by the (higher) Dehn twists along Langragian spheres.
On the algebraic geometry side,
Seidel-Thomas \cite{ST} studied the mirror counterpart of \cite{KS} (also in type A).
They showed that $\mathcal{D}(\qq{N})$ can be realized as a subcategory of the bounded derived category of
coherent sheaves of the mirror variety.
By their work, for any (acyclic) quiver $Q$, there is
the spherical twist group $\Br(\qq{N})$ on Calabi-Yau-$N$ category $\mathcal{D}(\qq{N})$.
Similarly one can defined the spherical twist groups for the quivers with (super)potential.

\subsection{Spherical twists on Calabi-Yau categories}\label{sec:ST-CY}
Fix an algebraically closed field $\k$.
An (algebraic) triangulated category $\hua{D}$ is Calabi-Yau-$N$
if for any objects $L,M$ in $\hua{D}$ there is a natural isomorphism
\begin{gather}\label{eq:serre}
    \mathfrak{S}:\Hom_{\hua{D}}^{\bullet}(L,M)
        \xrightarrow{\sim}\Hom_{\hua{D}}^{\bullet}(M,L)^\vee[N].
\end{gather}
An object $S$ in such a category is spherical if $\Hom^{\bullet}(S, S)=\k \oplus \k[-N]$
and it induces a spherical twist functor $\phi_S\in\Aut\D$,
defined by
\begin{equation}\label{eq:sphtwist+}
    \phi_S(X)=\Cone\left(S\otimes\Hom^\bullet(S,X)\to X\right).
\end{equation}

Consider a (graded) quiver with superpotential $(Q,W)$ of degree $N$,
in the sense of \cite{G,K8,VB,O}.
For instance, it satisfies at least the following properties:
\begin{itemize}
\item the degrees of arrows of $Q$ are in $\{0,\ldots,1-N\}$;
\item $W$ is the sum of homogenous cycles of degree $3-N$;
\item there is a (distinguish) loop of degree $2-N$ at each vertices.
\end{itemize}
We also require it is good,
i.e. it satisfies certain conditions, e.g. \cite[Section~6]{O} for details.
Denote by $$\Gamma\colon=\Gamma(Q,W)$$ the \emph{Ginzburg dg $\k$-algebra (of degree N)}
associated to $(Q,W)$ and $\D_{fd}(\Gamma)$ the finite-dimensional derived category of $\Gamma$.
Then $\D_{fd}(\Gamma)$ is a Calabi-Yau-$N$, Hom-finite, Krull-Schmidt, $\k$-linear triangulated category.
We also know that $\D_{fd}(\Gamma)$ admits a canonical heart $\zero$ generated
by the simple $\Gamma$-modules $\{S_i\mid i\in Q_0\}$.

\begin{definition}\label{def:ST}
The spherical twist group $\ST(Q,W)$ of a quiver with superpotential
$(Q,W)$ is the subgroup of $\Aut\D_{fd}(Q,W)$ generated by the spherical twists of
the simple $\Gamma$-modules.
\end{definition}

\begin{remark}
Our main motivation is to study the spherical twist groups,
which has a closed relation with proving the
topological properties (namely simply connectedness and contractibility) of
the corresponding stability space of $\mathcal{D}(\qq{N})$, in the sense of Bridgeland (\cite{B1,Q2,QW}).
Our previously works \cite{Q2,QQ,QZ2,QZ3} are attempts to understand the spherical twist groups
via other generalizations of braid groups.
\end{remark}

\subsection{Relations between different braid groups}\label{sec:pre}
Here we summarize results on the relations between various generalized braid groups.
Recall that we have the following generalizations:
\begin{itemize}
  \item Braid group $\Br$ in Definition~\ref{def:Br};
  \item Algebraic braid twist group $\AT$ in Definition~\ref{def:AT};
  \item Braid twist group $\BT$ in Definition~\ref{def:BT};
  \item Dehn twist group $\DT$ in Definition~\ref{def:DT};
  \item Spherical twist group $\ST$ in Definition~\ref{def:ST}.
\end{itemize}

\begin{theorem}
Let $(Q,W)$ be a quiver with potential.
We have the following.
\begin{itemize}
  \item[\cite{KS,ST}] If $(Q,W)=(A_n,0)$ and $N\geq2$ is an integer, then
  \begin{gather}\label{eq:ST}
    \Br_{n+1}=\Br(A_n)\cong\ST(Q,W).
  \end{gather}
  \item[\cite{BT}] If $(Q,W)=(Q,0)$ is a Dynkin quiver and $N=2$, then
  \eqref{eq:ST} holds.
  \item[\cite{IUU}] If $(Q,W)=(\widetilde{A_n},0)$ and $N=2$, then
  \eqref{eq:ST} holds.
  \item[\cite{QW}] If $(Q,W)=(Q,0)$ is a Dynkin quiver and $N\geq2$ is an integer, then
  \eqref{eq:ST} holds.
  \item[\cite{QQ}] If $(Q,W)$ is coming from a triangulation of
    an unpuncture marked surface $\surf$ and $N=3$, then
  \begin{gather}\label{eq:QQ}
    \BT(Q,W)\cong\ST(Q,W).
  \end{gather}
  \item[\cite{GM,QQ}] If $(Q,W)$ is mutation (in the sense of \cite{FZ}) equivalent to
    a Dynkin quiver $Q^*$ and $N=3$, then
  \begin{gather}\label{eq:GM}
    \AT(Q,W)\cong\Br(Q^*).
  \end{gather}
  \item [\cite{QZ3}] If $(Q,W)$ is coming from a triangulation of
    an unpuncture marked surface $\surf$ and $N=3$ and there is no double arrows in $Q$, then
  \begin{gather}\label{eq:QZ}
    \AT(Q,W)\cong\BT(Q,W).
  \end{gather}
  Moreover, when there are double arrows in $Q$,
  \cite{QZ3} also gives the generalization/modification of the definition of $\AT(Q,W)$
  so that the equality above still holds.
  This is equivalent to give the presentation of $\BT(Q,W)$.
  \item[\cite{KQ2}] If $(Q,W)$ is coming from a triangulation of
    an unpuncture marked surface $\surf$ and $N=3$, then
  \begin{gather}\label{eq:KQ}
    \DT(Q,W)\cong\BT(Q,W).
  \end{gather}
  If $(Q,W)$ is mutation equivalent to a type $D_n$ quiver (for $\surf$ is a once-punctured $n$-gon), then
  \[
    \DT(Q,W)\cong\Br(D_n).
  \]
\end{itemize}

\end{theorem}

\section{Cluster braid groups}\label{sec:4}
Now we introduce another generalization of braid groups for quivers with superpotential,
which is due to the joint work with Alastair King \cite{KQ2}.
\subsection{Exchange graphs and exchange groupoids}
Fix an integer $N\geq3$.
Consider a (good) quiver with superpotential $(Q,W)$ of degree $N$
with the associated Ginzburg dg algebra $\Gamma=\Gamma(Q,W)$.
We have the following associated categories (cf. \cite{K8,KQ})
with certain exchange (oriented) graph structures:
\begin{itemize}
\item The Calabi-Yau-$N$ category $\D_{fd}(\Gamma)$
with hearts (as vertices) and simple forward tiltings (as edges);
\item The perfect derived category $\per\Gamma$
with silting objects (as vertices) and forward mutations (as edges);
\item The higher ($N-1$) cluster category $\C(\Gamma)$
with cluster tilting objects (as vertices) and forward mutations (as edges).
\end{itemize}
Denote (the principal components of) these three exchange graphs by
$\EG(\Gamma)$, $\SEG(\Gamma)$ and $\CEG{N-1}{\Gamma}$ respectively.
We have the following result (cf. \cite[Section~3.1]{QQ2} and comments there).
\begin{gather}\label{eq:cong1}
    \EG(\Gamma)\cong\SEG(\Gamma),\\\label{eq:cong2}
    \EG(\Gamma)/\ST(Q,W)\cong\CEG{N-1}{\Gamma}.
\end{gather}

For convenience, we will work with Ext quivers of hearts 
(where the definition is more straightforward)
instead of original quivers in the quivers with superpotential setting.
They differ by grading shift and doubling (\cite[Definition~6.1 and Theorem~8.10]{KQ}).

\begin{definition}\label{def:extquiv}
Let $\h$ be a finite heart in a triangulated category $\D$.
The \emph{Ext quiver} $\hua{Q}_{\h}$ is the (positively) graded quiver
whose vertices are the simples of $\h$ and
whose degree $k$ arrows $S_i \to S_j$ correspond to a basis of
$\Hom_{\D}(S_i, S_j[k])$.
For a cluster $\mathbf{C}$, the associated quiver $\hua{Q}_\mathbf{C}$
is defined to the Ext quiver of any its lifted heart in $\EG(\Gamma)$.
This is well-defined since auto-equivalences preserve Ext quivers.
\end{definition}

The generators of fundamental groups of these graphs are essentially
squares, pentagons and $2N$-gons.
In the finite type case, the $2N$-gons are actually covered by squares and pentagons
(cf. \cite[Proof of Theorem~5.4]{Q2}).
Denote by $\tilt{\h}{\sharp}{S}$ the simple forward tilting of $\h$ w.r.t. a simple $S$
(where the torsion part in the torsion pair in $\h$ for the tilting is generated by $S$).
The inverse operation is the backward simple tilting, denoted by $\tilt{\h}{\flat}{S}$,
and we have $\h=\tilt{\left( \tilt{\h}{\sharp}{S} \right)}{\flat}{S[1]}$.
And inductively we define that,
\[
    \tilt{\h}{m\sharp }{S}
    ={  \Big( \tilt{\h}{ (m-1) \sharp}{S} \Big)  }^{\sharp}_{S[m-1]} \, 
\]
for $m\geq1$, and similarly we have $\tilt{\h}{m\flat }{S}$ for $m\geq1$.
For $m<0$, we also set $\tilt{\h}{m\sharp }{S}=\tilt{\h}{-m\flat }{S}$.
Note that we have (cf. \cite[Cororllary~8.4]{KQ})
    \begin{equation}\label{eq:root}
        \tilt{\h}{ \pm(N-1)\sharp}{S} = \phi_{S}^{\mp1} (\h ).
    \end{equation}
More details on (Happel-Reiten-Smal\o) tilting theory can be found in \cite[Section~3]{KQ}.
We have the following result.

\begin{figure}[ht]
\begin{tikzpicture}[arrow/.style={->,>=stealth},rotate=0,scale=1.35]
\draw (0,0) node (t1) {$\h$}
    (0,2) node (t2) {$\h_j$}
    (2,2) node (t3) {$\h_{ji}$}
    (4,2) node (t4) {$\tilt{\left(\h_j\right)}{m\sharp}{S_i}$} (6,2) node (tt4) {$\phi_{S_i}^{-1}(\h_j)$}
    (4,0) node (t5) {$\tilt{\h}{m\sharp}{S_i}$} (6,0) node (tt5) {$\phi_{S_i}^{-1}(\h)$}
    (2,0) node (t6) {$\h_i$};
\draw [->, font=\scriptsize]
    (t1)edge node[left]{$S_j$} (t2)
    (t2)edge node[above]{$S_i$} (t3)
    (tt5)edge node[right]{$T_j$} (tt4)
    (t6)edge[dashed] node[above]{} (t5)
    (t5)edge[dashed] node[above]{} (tt5)
    (t1)edge node[below]{$S_i$} (t6)
    (t3)edge[dashed] node[below]{} (t4)
    (t4)edge[dashed] node[below]{} (tt4);
\end{tikzpicture}

\begin{tikzpicture}[arrow/.style={->,>=stealth},rotate=0,scale=1.35]
\draw (0,0) node (t1) {$\h$}
    (0,2) node (t2) {$\h_j$}
    (2,2) node (t3) {$\h_{ji}$}
    (2,0) node (t6) {$\h_i$};
\draw [->, font=\scriptsize]
    (t1)edge node[left]{$S_j$} (t2)
    (t2)edge node[above]{$S_i$} (t3)
    (t1)edge node[below]{$S_i$} (t6)
    (t6)edge node[right]{$S_j$} (t3);
\end{tikzpicture}
\quad
\begin{tikzpicture}[scale=1.35,arrow/.style={->,>=stealth},rotate=0]
\draw (0,0) node (t1) {$\h$}
    (0,2) node (t2) {$\h_j$}
    (2,2) node (t3) {$\h_{ji}$}
    (2,0) node (t6) {$\h_i$}
    (2,1) node (t) {$\h_*$};
\draw [->, font=\scriptsize]
    (t1)edge node[left]{$S_j$} (t2)
    (t2)edge node[above]{$S_i$} (t3)
    (t1)edge node[below]{$S_i$} (t6)
    (t6)edge node[right]{$T_j$}(t)
    (t)edge node[right]{$S_j$}(t3);
\end{tikzpicture}
\caption{Square, pentagon and $2N$-gon}
\label{fig:hex}
\end{figure}
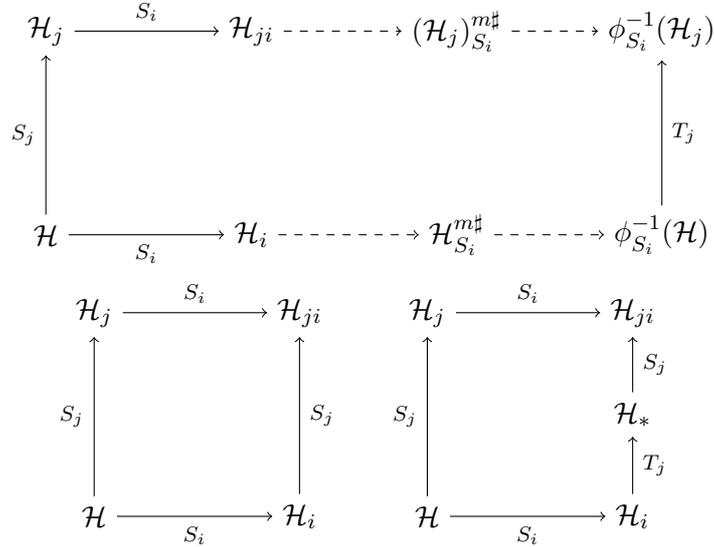

\begin{proposition}\cite{KQ2,Q2,QW}\label{pp:eg.h}
Let $\h$ be a heart in $\EG(\Gamma)$ with simples $S_i$ and $S_j$ satisfying
$\Ext^1(S_i,S_j)=0$ and $\h_i=\tilt{\h}{\sharp}{S_i},\h_j=\tilt{\h}{\sharp}{S_j}$.
We have the following.
\begin{itemize}
\item[(O).] $S_i$ is a simple in $\h_j$ and there is an (oriented) $2N$-gon in $\EG(\Gamma)$,
as shown in the upper picture of Figure~\ref{fig:hex},
where $$\h_{ji}\colon=\tilt{ (\h_j) }{\sharp}{S_i},\quad T_j\colon=\phi_{S_i}^{-1}(S_j).$$
\item[(I).]
    If further $\Ext^1(S_j, S_i)=0$, then $T_j=S_j$, $\tilt{ (\hua{H}_i) }{\sharp}{S_j}=\hua{H}_{ji}$
    and there is a square as shown in the lower left picture of Figure~\ref{fig:hex}.
\item[(II).]
    If further $\Ext^1(S_j, S_i)=\mathbf{k}$, then  
    there is a pentagon as shown in the lower right picture of Figure~\ref{fig:hex},
    where $$\hua{H}_{ij}\colon=\tilt{  (\hua{H}_*)  }{\sharp}{S_j},\quad 
        \hua{H}_*\colon=\tilt{ (\hua{H}_i) }{\sharp}{T_j}.$$
\end{itemize}
\end{proposition}

By the isomorphisms \eqref{eq:cong1}, \eqref{eq:cong2},
there are such squares, pentagons and $2N$-gons in $\SEG(\Gamma)$ as well (cf. \cite{KQ2}).
However we need to manually add these relations to $\CEG{N-1}{\Gamma}$.
Let $\mathrm{E}$ be an oriented graph.
Denote by $\hua{W}^+(\mathrm{E})$ the \emph{path category} of $\mathrm{E}$,
i.e. whose objects are the vertices of $\mathrm{E}$ and whose generating morphisms are
the (oriented) edges of $\mathrm{E}$.
Denote by $\hua{W}(\mathrm{E})$ the \emph{path groupoid} of $\mathrm{E}$,
i.e. the same presentation of $\hua{W}^+(\mathrm{E})$ but all the morphisms are invertible.

\begin{definition}[King-Qiu \cite{KQ2}]\label{def:gpd}
The \emph{exchange groupoid} $\eg(\Gamma)$ is defined
to be the quotient groupoid of the path groupoid $\hua{W}(\EG(\Gamma))$
by the square, pentagon and $2N$-gon relations as in Proposition~\ref{pp:eg.h}.
Similarly,
the \emph{cluster exchange groupoid} $\ceg_{N-1}(\Gamma)$
is the quotient groupoid of the path groupoid $\hua{W}(\CEG{N-1}{\Gamma})$
by the induced square, pentagon and $2N$-gon relations in Proposition~\ref{pp:eg.h} (via \eqref{eq:cong2}).
\end{definition}

\begin{definition}[King-Qiu \cite{KQ2}]\label{def:CBr}
Consider a quiver with superpotential $(Q,W)$ of degree $N$
with the associated Ginzburg dg algebra $\Gamma$.
Its cluster braid group $\CBr(Q,W)$
is the point group of the cluster exchange groupoid $\ceg_{N-1}(\Gamma)$, i.e.
$$\CBr(Q,W)=\CBr(\mathbf{C}_{\Gamma})\colon=\pi_1 (\ceg_{N-1}(\Gamma),\mathbf{C}_{\Gamma}),$$
where $\mathbf{C}_\Gamma$ is the canonical cluster tilting object induced by $\Gamma$.
\end{definition}

\begin{remark}[Generators]
Here are some more detailed description of the generators of the point group $\CBr(Q,W)$.
By formula \eqref{eq:root}, there is a length $N-1$ path
\[
    \h\xrightarrow{S} \tilt{\h}{\sharp}{S}\xrightarrow{S[1]}\cdots\to\tilt{\h}{(N-1)\sharp}{S}
    =\phi_S^{-1}(\h)
\]
in $\EG(\Gamma)$
and by \eqref{eq:cong1}, it becomes a $(N-1)$-loop $l_S$ at the corresponding vertex,
some cluster tilting object $\mathbf{C}$ in $\CEG{N-1}{\Gamma}$.
Moreover, each simple of $\h$ corresponds to an indecomposable summand $Y$ of $\mathbf{C}$,
which corresponds to a (forward) mutation.
In fact, the loop $l_S$ contains all cluster tilting objects, which are
completions of the almost complete cluster tilting object $\mathbf{C}-Y$.
In other words, this type of loops is indexed by almost complete cluster tilting objects.

Then the generators of a cluster braid group $\CBr(\mathbf{C})$ are the loops
indexed by almost complete cluster tilting objects, which are summand of $\mathbf{C}$.
Locally, they can be labelled by vertices of the Ext quiver 
$\hua{Q}_\mathbf{C}$ associated to $\mathbf{C}$. 
\end{remark}

\begin{remark}[Conjugation formula for Calabi-Yau-$N$]\label{rem:conj}
One of the key to generalize the result/construction in \cite{KQ2} to the CY-$N$ setting is
the following conjugation formula.
Given a forward mutation $\mathbf{C}\xrightarrow{x}\mathbf{C}'$ in $\ceg_{N-1}(\Gamma)$.
Note that (locally) we can identify the vertex sets of $Q_\mathbf{C}$ and $Q_{\mathbf{C}'}$
(say $\{i\}$)
and suppose that $x=\mu_j$ is w.r.t. vertex $j$ (i.e. $x=\mu_j$ at $\mathbf{C}$).
Denote by $\{t_i\}$ the local twists/generators of $\CBr(\mathbf{C})$
and by $\{t_i'\}$ the local twists/generators in $\CBr(\mathbf{C}')$.
The conjugation of $x$ in $\ceg_{N-1}(\Gamma)$ gives an isomorphism
\[\begin{array}{rcl}
    \conj_x\colon \CBr(\mathbf{C})&\to&\CBr(\mathbf{C}')\\
    t&\mapsto& x^{-1} t x.
\end{array}\]
such that
\begin{gather*}
\begin{array}{ccl}
    \conj_x(t_i)&=&\begin{cases}
    (t'_j)^{-1} t_i' t_j'& \text{if there are arrows of degree $1$ from $i$ to $j$ in $\hua{Q}_\mathbf{C}$,}\\
    t_i' & \text{otherwise,}
    \end{cases}\end{array}
\end{gather*}
The proof is basically using $2N$-gon relation as in \cite{KQ2}.
The non-trivial calculation is as follows.
Suppose that there are arrows of degree $1$ from $i$ to $j$ in $\hua{Q}_\mathbf{C}$.
Consider a lift/heart $\h$ of $\mathbf{C}$ in $\EG(\Gamma)$
and denote the corresponding simples by $S_i$ and $S_j$.
So the lift of the edge $\mathbf{C}\xrightarrow{x}\mathbf{C}'$ in $\EG(\Gamma)$ is the tilting
$$\h\xrightarrow{S_j}\tilt{\h}{\sharp}{S_j}=\colon\h_j.$$
Then in the following tilting sequence
\[
    \phi_{S_j}(\h_j)=\tilt{\h}{(2-N)\sharp}{S_j} \to\cdots\xrightarrow{S_j[-1]}
        \h\xrightarrow{S_j} \h_j,
\]
all of the hearts except $\h_j$ have $S_i$ as the $i$-th simple
(while $\h_j$ has the $i$-th simple $S'_i=\phi^{-1}_{S_j}(S_i)$.)
Then we have the following sub-graph of $\EG(\Gamma)$,
\begin{figure}[ht]
\begin{tikzpicture}[arrow/.style={->,>=stealth},xscale=1.4]
\draw (0,0) node (t1) {$\h$}
    (0,2) node (t2) {$\h_j$}
    (2,2) node (t3) {}
    (4,2) node (t4) {} (6,2) node (tt4) {$\phi_{S_i}^{-1}(\h_j)$}
    (4,0) node (t5) {$\tilt{\h}{m\sharp}{S_i}$} (6,0) node (tt5) {$\phi_{S_i}^{-1}(\h)$}
    (2,0) node (t6) {$\h_i$}

    (0,-2) node (a1) {$\tilt{\h}{\flat}{S_j}$}
    (4,-2) node (a5) {} (6,-2) node (at5) {$\phi_{S_i}^{-1}(\tilt{\h}{\flat}{S_j})$}
    (2,-2) node (a6) {}
    (0,-4) node (b1) {$\tilt{\h}{k\flat}{S_j}$}
    (4,-4) node (b5) {} (6,-4) node (bt5) {$\phi_{S_i}^{-1}(\tilt{\h}{k\flat}{S_j})$}
    (2,-4) node (b6) {}
    (0,-6) node (c1) {$\phi_{S_j}(\h)$}
    (4,-6) node (c5) {} (6,-6) node (ct5) {$\phi_{S_i}^{-1}(\phi_{S_j}(\h))$}
    (2,-6) node (c6) {};
\draw [->, font=\scriptsize]
    (t1)edge[blue] node[left]{$S_j$} (t2)
    (tt5)edge[blue] node[right]{$\phi_{S_i}^{-1}(S_j)$} (tt4)

    (t6)edge[red,dashed] node[above]{} (t5)
    (t5)edge[red,dashed] node[above]{} (tt5)
    (t1)edge[red] node[above]{$S_i$} (t6)

    (a6)edge[dashed] node[above]{} (a5)
    (a5)edge[dashed] node[above]{} (at5)
    (a1)edge node[above]{$S_i$} (a6)

    (b6)edge[dashed] node[above]{} (b5)
    (b5)edge[dashed] node[above]{} (bt5)
    (b1)edge node[above]{$S_i$} (b6)

    (c6)edge[JungleGreen,dashed] node[above]{} (c5)
    (c5)edge[JungleGreen,dashed] node[above]{} (ct5)
    (c1)edge[JungleGreen] node[above]{$S_i$} (c6)

    (a1)edge[JungleGreen] node[left]{$S_j[-1]$}(t1)
        (at5)edge[JungleGreen] node[right]{$\phi_{S_i}^{-1}(S_j[-1])$}(tt5)

    (b1)edge[JungleGreen,dashed](a1)(c1)edge[JungleGreen,dashed](b1)
    (bt5)edge[JungleGreen,dashed](at5)(ct5)edge[JungleGreen,dashed](bt5);
\end{tikzpicture}
\caption{$N-2$ distinguish $2N$-gons}
\label{fig:cool}
\end{figure}
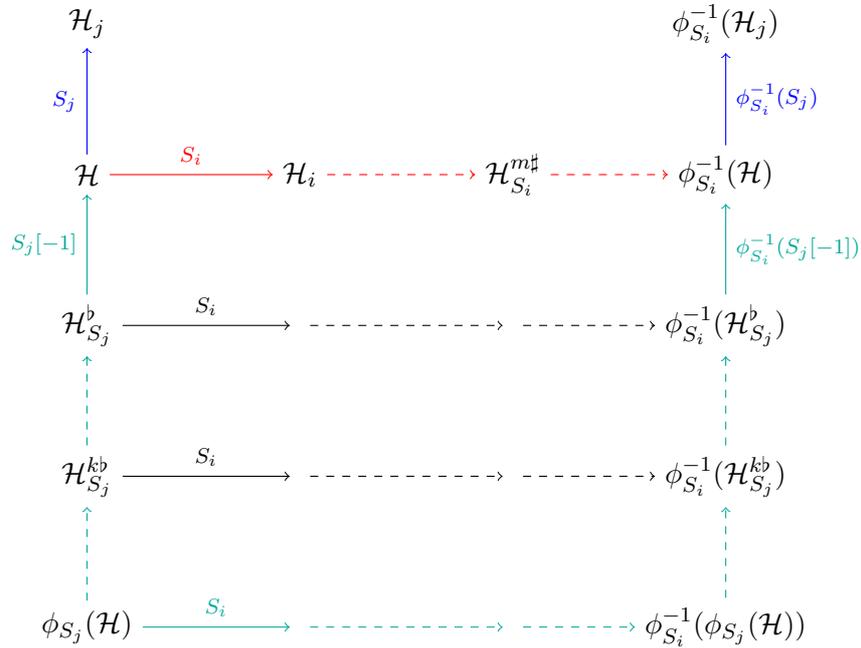
which consists of $(N-2)$ distinguish $2N$-gons.
Note that in Figure~\ref{fig:cool}, the right column is exactly
the twist $\phi_{S_i}^{-1}$ of the left column;
and except the first row, such a twist is realized by $N-1$ simple forward tiltings
by formula \eqref{eq:root}.

The lift of $x^{-1} t_i x$ in $\EG(\Gamma)$
is the path in Figure~\ref{fig:cool} consisting of blue and red edges;
the lift of $(t'_j)^{-1} t_i' t_j'$ in $\EG(\Gamma)$
is the path in Figure~\ref{fig:cool} consisting of blue and green edges.
As they differ by $(N-2)$ $2N$-gons, we have
$x^{-1} t_i x=(t'_j)^{-1} t_i' t_j'$ as required.
\end{remark}

\subsection{Two examples: $A_2$ and $A_1\times A_1$}
Fix $N=3$ in this subsection.
Notice that our definition of cluster exchange graph is different from the usual.
More precisely, our $\CEG{2}{\Gamma}$ can be obtained from the usual
cluster exchange graph by replacing each unoriented edge with a 2-cycle.
This idea had appeared in \cite[Section~9]{KQ}.

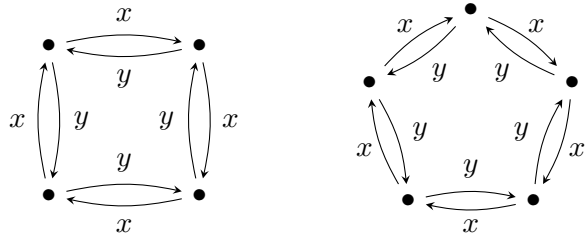
\begin{figure}[h]\centering
\newcommand{\vsource}{\otimes}
\newcommand{\vsink}{\odot}
\newcommand{\vertx}{\bullet}
\begin{tikzpicture}[scale=.7, rotate=45,
  arrow/.style={->,>=stealth}]
\foreach \j in {1,2,3,4}
   {\draw (90*\j:2cm) node (t\j) {$\bullet$};
   \draw (90*\j+45:2cm) node {$x$};\draw (90*\j+45:.8cm) node {$y$};}
\foreach \a/\b in {1/2,2/3,3/4,4/1}{
  \draw (t\a) edge[arrow,bend left=12]  (t\b);
  \draw (t\b) edge[arrow,bend left=12]  (t\a);}
\end{tikzpicture}
\qquad\quad
\begin{tikzpicture}[scale=.7, rotate=90,
  arrow/.style={->,>=stealth}]
\foreach \j in {1,2,3,4,5}
   {\draw (72*\j:2cm) node (t\j) {$\bullet$};
   \draw (72*\j+36:2.1cm) node {$x$};\draw (72*\j+36:1cm) node {$y$};}
\foreach \a/\b in {1/2,2/3,3/4,4/5,5/1}{
  \draw (t\a) edge[arrow,bend left=12]  (t\b);
  \draw (t\b) edge[arrow,bend left=12]  (t\a);}
\end{tikzpicture}
\caption{$\CEG{2}{\Gamma(Q)}$ for a quiver $Q$ of type $A_1\times A_1$ and $A_2$}
\label{fig:4}
\end{figure}

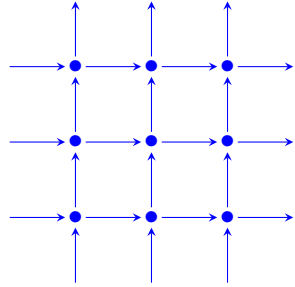
\begin{figure}[h]\centering
\newcommand{\vsource}{\otimes}
\newcommand{\vsink}{\odot}
\newcommand{\vertx}{\bullet}
\begin{tikzpicture}[arrow/.style={->,>=stealth}]
\foreach \i in {2,3,4}{
\foreach \j in {1,2,3,4,5}
   {\draw[white] (\i,\j) node (t\j) { };}
\foreach \a/\b in {1/2,2/3,3/4,4/5}{
  \draw[blue] (t\a) edge[arrow] (t\b);
}}
\foreach \i in {2,3,4}{
\foreach \j in {1,2,3,4,5}
   {\draw[white] (\j,\i) node (t\j) { };}
\foreach \a/\b in {1/2,2/3,3/4,4/5}{
  \draw[blue] (t\a) edge[arrow] (t\b);
}}

\foreach \i in {2,3,4}{
\foreach \j in {2,3,4}
   {\draw[blue] (\j,\i) node (t\j) {$\bullet$};}}
\end{tikzpicture}
\caption{$\EG(\Gamma(A_1\times A_1))$}
\label{fig:cover1}
\end{figure}

\begin{example}
Consider the quiver $Q=A_1\times A_1$ (and $W=0$).
Then $\CEG{2}{\Gamma(A_1\times A_1)}$ is shown in the left picture of Figure~\ref{fig:4},
where there should be four square faces attached,
that correspond to the relations $$x^2=y^2.$$
Then we have
\[
    \CBr(A_1\times A_1)=\mathbb{Z}^2\cong\ST(A_1\times A_1)
\]
and its universal cover is $\EG(\Gamma(A_1\times A_1))$, shown in Figure~\ref{fig:cover1}.
\end{example}

\begin{figure}\centering
\begin{tikzpicture}[scale=0.6,
 arrow/.style={->,>=stealth,thick,blue}, 
 border/.style={Periwinkle,dotted,thick},
 c-vrtx/.style={blue},
 c-arc/.style={Emerald}]
\newcommand{\vrtx}{\bullet}
\coordinate (O) at (0,0);
\coordinate (S1) at (0,6) ;
\draw [border] (O) circle (6cm);
\draw[c-arc,thick] (S1) \foreach \j in {1,...,3}
    {arc(360/3-\j*360/3+180:360-\j*360/3:10.3923cm)}--cycle;
\draw[c-arc,semithick] (S1) \foreach \j in {1,...,6}
    {arc(360/6-\j*360/6+180:360-\j*360/6:3.4641cm)}--cycle;
\draw[c-arc] (S1) \foreach \j in {1,...,12}
    {arc(360/12-\j*360/12+180:360-\j*360/12:1.6077cm)}--cycle;
\foreach \j in {1,...,3}
{\draw[c-vrtx] (-90+120*\j:0.7cm) node (v\j) {$\vrtx$};
 \draw[c-vrtx] (-210+120*\j:0.7cm) node (w\j) {$\vrtx$};
 \draw[c-vrtx] (-90+120*\j:2.2cm) node (a\j) {$\vrtx$};
 \draw[c-vrtx] (-90+15+120*\j:3cm) node (b\j) {$\vrtx$};
 \draw[c-vrtx] (-90-15+120*\j:3cm) node (c\j) {$\vrtx$};}
\foreach \j in {1,...,3}
{\draw [arrow] (v\j) edge[bend left] (w\j);
 \draw [arrow] (a\j) edge[bend left] (b\j);
 \draw [arrow] (b\j) edge[bend left] (c\j);
 \draw [arrow] (c\j) edge[bend left] (a\j);}
\foreach \j in {1,...,3}
{\draw[c-vrtx] (60*\j*2-1-60:3.9cm) node (x1\j) {$\vrtx$};
 \draw[c-vrtx] (60*\j*2+1-120:3.9cm) node (x2\j) {$\vrtx$};}
\foreach \j in {1,...,3}
{\draw [arrow]  (v\j) edge[bend left] (a\j);
 \draw [arrow]  (a\j) edge[bend left] (v\j);
 \draw [arrow]  (b\j) edge[bend left] (x1\j);
 \draw [arrow]  (c\j) edge[bend left] (x2\j);
 \draw [arrow]  (x1\j) edge[bend left] (b\j);
 \draw [arrow]  (x2\j) edge[bend left] (c\j);}
\end{tikzpicture}
\quad
\begin{tikzpicture}[scale=0.4,
  egarrow/.style={->,>=stealth,thick,blue}, 
  vline/.style={dashed,gray},
  c-vrtx/.style={blue}]
\newcommand{\vrtx}{\bullet}
\draw[fill=gray!7,dotted]
    (-7, 2.5) to (5, 2.5) to (3, -3) to (-9,-3) -- cycle;
\draw[c-vrtx] (-4,0) node (A) {$\vrtx$};
\draw[c-vrtx] (-4,10) node (Ax) {$\vrtx$};
\draw[c-vrtx] (-4,4) node (Axx) {$\vrtx$};
\draw[c-vrtx] (1,1) node (C) {$\vrtx$};
\draw[c-vrtx] (1,6) node (Cx) {$\vrtx$};
\draw[c-vrtx] (0,-2) node (B) {$\vrtx$};
\draw[c-vrtx] (0,8) node (Bx) {$\vrtx$};
\draw[c-vrtx] (-6,0) node (D) {$\vrtx$};
\draw[c-vrtx] (-6,7) node (Dx) {$\vrtx$};
\draw [vline] (A) edge (Axx); \draw [vline] (Axx) edge (Ax);
\draw [vline] (C) edge (Cx); \draw [vline] (Cx) edge (1,10);
\draw [vline] (B) edge (Bx); \draw [vline] (Bx) edge (0,10);
\draw [vline] (D) edge (Dx); \draw [vline] (Dx) edge (-6,10);
\draw [egarrow] (A) edge[bend left] (D);
\draw [egarrow] (D) edge[bend left] (A);
\draw [egarrow] (A) edge[bend left] (C);
\draw [egarrow] (C) edge[bend left] (B);
\draw [egarrow] (B) edge[bend left] (A);
\draw [egarrow] (Axx) edge (Dx);
\draw [egarrow] (Axx) edge (Cx);
\draw [egarrow] (Cx) edge (Bx);
\draw [egarrow] (Bx) edge (Ax);
\draw [egarrow] (Dx) edge (Ax);
\end{tikzpicture}
\caption{The quotient graph $\EG(\Gamma(A_2))/[1]$
and the lifting (local demonstration)}\label{fig:cover2}
\end{figure}
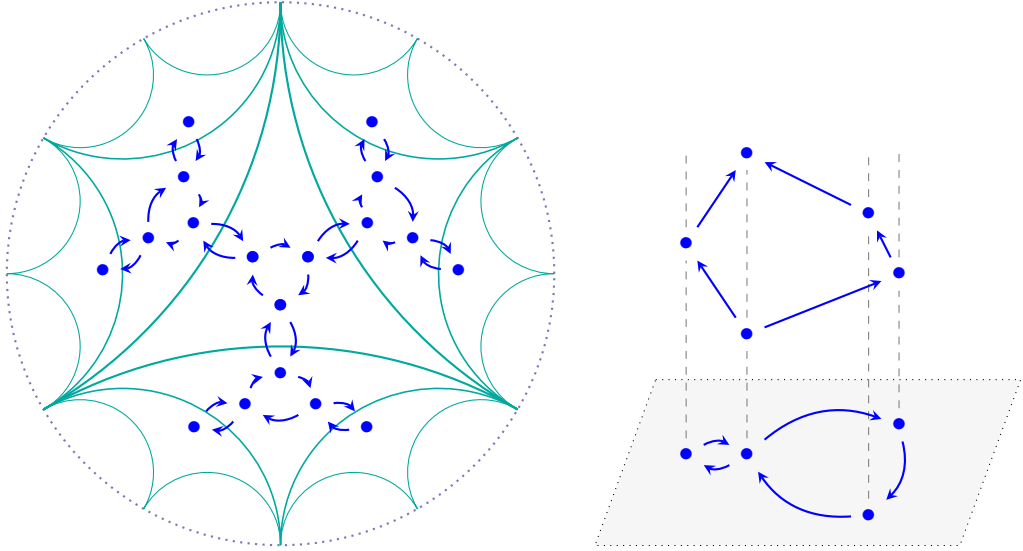

\begin{example}\cite[Section~10]{KQ}
Consider the quiver $Q=A_2$ (and $W=0$).
Then $\ceg_2(\Gamma(A_1\times A_1))$ is shown in the right picture of Figure~\ref{fig:4},
where there should be five pentagons faces attached.
that correspond to the relations $$x^2=y^3.$$
Then we have
\[
    \CBr(A_2)=\Br_3\cong\ST(A_2)
\]
and its universal cover is $\eg(\Gamma(A_2))$.
The quotient graph $\EG(\Gamma(A_2))/\mathbb{Z}[1]$ is shown
in the left picture of Figure~\ref{fig:cover2}
and its $\mathbb{Z}$-covering $\EG(\Gamma(A_2))$ can be constructed via lifting
shown in the right picture of Figure~\ref{fig:cover2}.
\end{example}

\subsection{Dynkin case}
In the Dynkin case, the phenomenon above also holds.
Namely, we have the following.

\begin{theorem}[Qiu-Woolf \cite{QW}]
Let $N\geq3$ be an integer and $(Q,W)$ be a quiver with superpotential
of Dynkin type (in the sense that $\Gamma(Q,W)$ is Morita equivalent to
the Calabi-Yau-$N$ completion of a Dynkin quiver $Q^*$).
Then we have
\begin{gather}
    \Br(Q^*)\cong\CBr(Q,W)\cong\ST(Q,W)(\subset\Aut\D_{fd}(\Gamma(Q,W)))
\end{gather}
and $\eg(\Gamma)$ is the universal cover of $\ceg_{N-1}(\Gamma)$.
Here the first isomorphism can be constructed inductively by choosing
a mutation sequence from $(Q,W)$ to $(Q^*,0)$ using the conjugation formula in Remark~\ref{rem:conj}.
\end{theorem}

\subsection{Decorated marked surface case}
Let $\surf$ be an unpunctured marked surface, $N=3$
and $(Q_\TT,W_\TT)$ be the quiver with potential associated to some triangulation $\TT$ of $\surfo$.

\begin{theorem}[King-Qiu \cite{KQ2}]
\begin{gather}
    \CBr(Q_\TT,W_\TT)\cong\BT(Q_\TT,W_\TT) \; (\cong\ST(Q_\TT,W_\TT) \text{ by \eqref{eq:QQ}}).
\end{gather}
and $\eg(\Gamma)$ is the universal cover of $\ceg_{2}(\Gamma)$.
As a consequence, the corresponding space of stability conditions $\Stap\D_{fd}(\Gamma))$ is simply connected.
\end{theorem}

\subsection{Conjectures in the general case}
In general, we expect the following:
\begin{conjecture}
For any quiver with superpotential $(Q,W)$ (of degree $N$),
\[
    \CBr(Q,W)\cong\ST(Q,W)
\]
and $\eg(\Gamma)$ is the universal cover of $\ceg_{N-1}(\Gamma)$.
\end{conjecture}
This conjecture is closely related to the conjectures that
the corresponding space $\Stap\D_{fd}(\Gamma)$ of stability conditions is
simply connected (and contractible) as appeared above in the surface case.



\end{document}